\documentclass{article}

\usepackage{amssymb,amsmath,amstext,amscd,amsthm,latexsym}

\numberwithin{equation}{section}

\theoremstyle{plain}
\newtheorem{theorem}{Theorem}

\newtheorem{remark}{Remark}
\newtheorem{example}{Example}

\begin{document}

\title{Infinite products of absolute zeta functions}
 
\author{Nobushige Kurokawa\footnote{Department of Mathematics, Tokyo Institute of Technology} \and Hidekazu Tanaka\footnote{6-15-11-202 Otsuka, Bunkyo-ku, Tokyo}}

\date{June 17, 2021}



\maketitle

\begin{abstract}
We study some infinite products of absolute zeta functions. Especially, we consider the convergence and the rationality of them.
\end{abstract}

\section{Introduction} It is well--known that some infinite products of zeta functions have interesting properties. For example the infinite product ${\displaystyle \prod_{n=1}^{\infty}\zeta(s+n)}$ has an analytic continuation to all $s \in {\mathbb C}$ as a meromorphic function, where $\zeta(s) = \zeta_{{\mathbb Z}}(s)$ is the Riemann zeta function: see Cohen--Lenstra \cite{CL} and Manin \cite[\S 3.5]{M}. We notice that 
\[
\prod_{n=1}^{\infty} \zeta(s+n) = \sum_{A} |{\rm Aut}(A)|^{-1} |A|^{-s},
\] 
where $A$ runs over all isomorphism classes of finite abelian groups. 
\par In this paper we study infinite products of absolute zeta functions (zeta functions over ${\mathbb F}_{1}$): see Soul{\'e} \cite{S}, Connes--Consani \cite{CC1} and Kurokawa--Ochiai \cite{KO}.
\par The simplest example is
\[
\prod_{n=1}^{\infty} \zeta_{{\mathbb F}_{1}}(s+n) = \prod_{n=1}^{\infty} \frac{1}{s+n},
\]
which diverges unfortunately.
\par We present some concrete examples of infinite products of absolute zeta functions such that
\[
\prod_{n=1}^{\infty} \zeta_{{\rm GL(3)}/{\mathbb F}_{1}}(s+n) = \zeta_{{\rm SL(3)}/{\mathbb F}_{1}}(s),
\]
where
\[
\zeta_{{\rm GL(3)}/{\mathbb F}_{1}}(s) = \frac{(s-8)(s-7)(s-3)}{(s-9)(s-5)(s-4)}
\]
and
\[
\zeta_{{\rm SL(3)}/{\mathbb F}_{1}}(s) = \frac{(s-6)(s-5)}{(s-8)(s-3)}.
\]
\par According to Kurokawa-Ochiai \cite{KO}, for a suitable scheme $X$ the absolute zeta function $\zeta_{X/{\mathbb F}_{1}}(s)$ is defined as 
\[
\zeta_{X/{\mathbb F}_{1}}(s)=\exp\biggl( \frac{\partial}{\partial w} Z_{X/{\mathbb F}_{1}}(w,s) \biggl|_{w=0} \biggl)
\]
with 
\[
Z_{X/{\mathbb F}_{1}}(w,s)=\frac{1}{\Gamma(w)} \int_{1}^{\infty} |X({\mathbb F}_{x})| x^{-s-1} (\log x)^{w-1} dx.
\] 
\begin{theorem} Let $r \geq 2$ be an integer. Then we have
\[
\zeta_{{\rm SL}(r)/{\mathbb F}_{1}}(s) = \prod_{n=1}^{\infty} \zeta_{{\rm GL}(r)/{\mathbb F}_{1}}(s+n).
\]
\end{theorem}

\begin{example} ($r=2$)
\begin{align*}
\zeta_{{\rm SL}(2)/{\mathbb F}_{1}}(s) &= \frac{s-1}{s-3},\\
\zeta_{{\rm GL}(2)/{\mathbb F}_{1}}(s) &= \frac{(s-2)(s-3)}{(s-1)(s-4)}
\end{align*}
and
\[
\zeta_{{\rm SL}(2)/{\mathbb F}_{1}}(s) = \prod_{n=1}^{\infty} \zeta_{{\rm GL}(2)/{\mathbb F}_{1}}(s+n).
\]
\end{example}

\begin{remark} ($r=1$)
\begin{align*}
\zeta_{{\rm SL}(1)/{\mathbb F}_{1}}(s) &= \frac{1}{s},\\
\zeta_{{\rm GL}(1)/{\mathbb F}_{1}}(s) &= \frac{s}{s-1}
\end{align*}
and
\[
\zeta_{{\rm SL}(1)/{\mathbb F}_{1}}(s) \neq \prod_{n=1}^{\infty} \zeta_{{\rm GL}(1)/{\mathbb F}_{1}}(s+n).
\]
\end{remark}

Let $N$ be a positive integer. For a zeta function $Z(s)$ we define a finite product $Z_{N}^{K}(s)$ and an infinite product $Z_{N}^{\infty}(s)$ as
\[
Z_{N}^{K}(s)=\prod_{k=0}^{K} Z(s+kN)
\]
and
\[
Z_{N}^{\infty}(s)=\prod_{k=0}^{\infty} Z(s+kN)
\]
respectively. Let 
\[
Z(s) = \zeta_{{\rm GL}(2)/{\mathbb F}_{1}}(s) = \frac{(s-2)(s-3)}{(s-1)(s-4)}
\]
in Theorems 2--5.

\begin{theorem}[explicit formula]
For $K=0,1,2,...$ we obtain the following results.
\par (1)
\[
Z_{1}^{K}(s) = \prod_{k=0}^{K}Z(s+k)=\frac{(s-2)(s+K-3)}{(s-4)(s+K-1)}.
\]
\par (2)
\begin{equation}
\label{1.1}
Z_{1}^{\infty}(s) = \frac{s-2}{s-4}.
\end{equation}
\par (3)
\[
Z_{2}^{K}(s)=\prod_{k=0}^{\infty}Z(s+2k)=\frac{(s-3)(s+2K-2)}{(s-4)(s+2K-1)}.
\]
\par (4)
\begin{equation}
\label{1.2}
Z_{2}^{\infty}(s)=\frac{s-3}{s-4}.
\end{equation}
\end{theorem}

\begin{remark} Let $q$ be a prime power. The congruence zeta function $\zeta_{X/{\mathbb F}_{q}}(s)$ associated with $X$ over ${\mathbb F}_{q}$ is defined by
\begin{align*}
\zeta_{X/{\mathbb F}_{q}}(s) = \exp \biggl( \sum_{m=1}^{\infty} \frac{|X({\mathbb F}_{q^m})|}{m}q^{-ms} \biggl).
\end{align*}
The Hasse zeta function $\zeta_{X/{\mathbb Z}}(s)$ associated with $X$ over ${\mathbb Z}$ is defined by
\[
\zeta_{X/{\mathbb Z}}(s) = \prod_{p:primes} \zeta_{X/{\mathbb F}_{p}}(s).
\]
In this case it is not difficult to show that $\zeta_{X/{\mathbb Z}}(s)$ has a functional equation
of the type
\[
\Gamma_{X/{\mathbb Z}}(s)\zeta_{X/{\mathbb Z}}(s) = \biggl( \Gamma_{X/{\mathbb Z}}(d-s)\zeta_{X/{\mathbb Z}}(d-s) \biggl)^{(-1)^{n_{1} + \cdots + n_{r}}}
\]
where $d=\sum_{j=1}^{r}\frac{n_{j}(3n_{j}-1)}{2}+1$ and $\Gamma_{X/{\mathbb Z}}(s)$ is expressed in terms of a product/quotient of the gamma function. We call $\Gamma_{X/{\mathbb Z}}(s)$ the gamma factor for $\zeta_{X/{\mathbb Z}}(s)$. In \cite{T}, We study the rationality of gamma factors associated to certain Hasse zeta functions. We refer to Manin \cite{M} and Connes-Consani \cite{CC2} for gamma factors of Hasse zeta functions. Especially, we have
\[
\Gamma_{{\rm GL}(2)/{\mathbb Z}}(s) = \frac{s-3}{s-4}.
\]
The right side is the same formula as the right side of (\ref{1.2}).
\end{remark}

\begin{theorem}[functional equation] For $K=0,1,2,...$ we obtain the following results.
\par (1)
\[
Z_{1}^{K}(5-K-s)=Z_{1}^{K}(s).
\]
\par (2)
\[
Z_{1}^{\infty}(6-s)=Z_{1}^{\infty}(s)^{-1}.
\]
\par (3)
\[
Z_{2}^{K}(5-2K-s)=Z_{2}^{K}(s).
\]
\par (4)
\[
Z_{2}^{\infty}(7-s)=Z_{2}^{\infty}(s)^{-1}.
\]
\end{theorem}

\begin{theorem}[convergence] $Z_{N}^{\infty}(s)$ converges and
\[
Z_{N}^{\infty}(s) = \frac{\Gamma(\frac{s-1}{N})\Gamma(\frac{s-4}{N})}{\Gamma(\frac{s-2}{N})\Gamma(\frac{s-3}{N})}.
\]
\end{theorem}

\begin{theorem}[rationality] $Z_{N}^{\infty}(s)$ is a rational function if and only if $N=1,2$.
\end{theorem}

Next, let 
\[
Z(s) = \frac{(s-c)(s-d)}{(s-a)(s-b)}.
\] in Theorems 6--8. 

\begin{theorem} The following (1) and (2) are equivalent.
\par (1) $Z_{N}^{\infty}(s)$ converges.
\par (2) $a+b=c+d$.
\end{theorem}

\begin{theorem}
The following (1), (2) and (3) are equivalent.
\par (1) $Z_{N}^{\infty}(s)$ converges.
\par (2) [functional equation] $Z(a+b-s)=Z(s)$.
\par (3) [absolute automorphy] $f(\frac{1}{x})=x^{-(a+b)}f(x)$. 
\end{theorem}

\begin{theorem}[rationality] $Z_{N}^{\infty}(s)$ is a rational function if and only if $\zeta_{N}^{a}+\zeta_{N}^{b}=\zeta_{N}^{c}+\zeta_{N}^{d}$ with $\zeta_{N}=\exp(\frac{2 \pi i}{N})$.
\end{theorem}

Finally, let 
\[
Z(s) = \zeta_{{\rm G}_{m}^{r}/{\rm F}_{1}}(s)=\prod_{\ell=0}^{r}(s-\ell)^{(-1)^{r-\ell+1}
{\scriptsize \left(
\begin{array}{c}
r\\
\ell
\end{array}
\right)}
}.
\] in Theorems 9--10. 

\begin{theorem} The following (1) and (2) are equivalent.
\par (1) $Z_{N}^{\infty}(s)$ converges.
\par (2) $r \geq 2$.
\end{theorem}

\begin{theorem}[rationality] $Z_{N}^{\infty}(s)$ is a rational function if and only if $N=1$.
\end{theorem}

\section{Proof of Theorem 1.}
\begin{proof}[Proof of Theorem 1] For $r \geq 2$ let
\[
f(x) = x^{r^2 - 1} (1-x^{-2}) \cdots (1-x^{-r}) = |{\rm SL}(r,{\mathbb F}_{x})|
\]
and
\[
g(x) = x^{r^{2}} (1-x^{-1}) \cdots (1-x^{-r}) = |{\rm GL}(r,{\mathbb F}_{x})|.
\]
Notice that $g(x)=(x-1)f(x)$. Then using $f(1)=g(1)=0$ we have (see \cite{KT1})
\[
\zeta_{{\rm SL}(r)/{\mathbb F}_{1}}(s) = \exp \biggl( \int_{1}^{\infty} \frac{f(x) x^{-s-1}}{\log x} dx \biggl)
\]
and
\[
\zeta_{{\rm GL}(r)/{\mathbb F}_{1}}(s) = \exp \biggl( \int_{1}^{\infty} \frac{g(x) x^{-s-1}}{\log x} dx \biggl).
\]
From these expressions we have
\begin{align*}
\frac{\zeta_{{\rm SL}(r)/{\mathbb F}_{1}}(s+n-1)}{\zeta_{{\rm SL}(r)/{\mathbb F}_{1}}(s+n)} &= \exp \biggl( \int_{1}^{\infty} \frac{f(x) (x^{-s-n}-x^{-s-n-1})}{\log x} dx \biggl)\\
&= \exp \biggl( \int_{1}^{\infty} \frac{f(x) (x-1)x^{-s-n-1}}{\log x} dx \biggl)\\
&= \exp \biggl( \int_{1}^{\infty} \frac{g(x) x^{-s-n-1}}{\log x} dx \biggl)\\
&= \zeta_{{\rm GL}(r)/{\mathbb F}_{1}}(s+n)
\end{align*}
for $n \geq 1$. Hence we get
\begin{align*}
\prod_{n=1}^{N} \zeta_{{\rm GL}(r)/{\mathbb F}_{1}}(s+n) &= \prod_{n=1}^{N} \frac{\zeta_{{\rm SL}(r)/{\mathbb F}_{1}}(s+n-1)}{\zeta_{{\rm SL}(r)/{\mathbb F}_{1}}(s+n)}\\
&= \frac{\zeta_{{\rm SL}(r)/{\mathbb F}_{1}}(s)}{\zeta_{{\rm SL}(r)/{\mathbb F}_{1}}(s+N)}.
\end{align*}
Thus it is sufficient to show that
\[
\lim_{N \to \infty} \zeta_{{\rm SL}(r)/{\mathbb F}_{1}}(s+N) = 1
\]
for $r \geq 2$. Let $f(x)=\sum_{k} a(k) x^{k} \in {\mathbb Z}[x]$ then
\begin{align*}
\zeta_{{\rm SL}(r)/{\mathbb F}_{1}}(s) &= \prod_{k} (s-k)^{-a(k)} \\
&= s^{-\sum_{k} a(k)} \prod_{k} (1-\frac{k}{s})^{-a(k)}.
\end{align*}
By using $f(1)=0=\sum_{k}a(k)$ we obtain
\[
\zeta_{{\rm SL}(r)/{\mathbb F}_{1}}(s)=\prod_{k} (1-\frac{k}{s})^{-a(k)}.
\]
Hence the expression
\[
\zeta_{{\rm SL}(r)/{\mathbb F}_{1}}(s+N)=\prod_{k} (1-\frac{k}{s+N})^{-a(k)}
\]
implies
\[
\lim_{N \to \infty} \zeta_{{\rm SL}(r)/{\mathbb F}_{1}}(s+N) = 1.
\]
Hence, we obtain Theorem 1.
\end{proof}

\section{Proof of Theorem 2--5.}
\begin{proof}[Proof of Theorem 2] 
\par (1) We prove (1) by the induction on $K=0,1,2,...$. Put $K=0$. Then we have
\[
Z_{1}^{0} = Z(s) = \frac{(s-2)(s-3)}{(s-4)(s-1)}.
\]
Now we assume 
\[
Z_{1}^{K}(s) = \frac{(s-2)(s+K-3)}{(s-4)(s+K-1)}.
\]
Then we have 
\begin{align*}
Z_{1}^{K+1}(s) &= Z_{1}^{K}(s)Z(s+K+1)\\
&=\frac{(s-2)(s+K-3)}{(s-4)(s+K-1)} \cdot \frac{(s+K-1)(s+K-2)}{(s+K)(s+K-3)}\\
&=\frac{(s-2)(s+K-2)}{(s-4)(s+K)}.
\end{align*}
Hence, we obtain (1).
\par (2) Let $K \to \infty$ in (1). Then we have
\[
\lim_{K \to \infty} \frac{s+K-3}{s+K-1}=1.
\]
Thus we obtain
\[
Z_{1}^{\infty}(s) = \frac{s-2}{s-4}.
\]
\par (3) We prove (3) by the induction on $K=0,1,2,...$. Put $K=0$. Then we have
\[
Z_{2}^{0}(s) = Z(s) =\frac{(s-3)(s-2)}{(s-4)(s-1)}.
\]
Now we assume 
\[
Z_{2}^{K}(s) = \frac{(s-3)(s+2K-2)}{(s-4)(s+2K-1)}.
\]
Then we have
\begin{align*}
Z_{2}^{K+1}(s) &= Z_{2}^{K}(s)Z(s+2(K+1))\\
&=\frac{(s-3)(s+2K-2)}{(s-4)(s+2K-1)} \cdot \frac{(s+2K)(s+2K-1)}{(s+2K+1)(s+2K-2)}\\
&=\frac{(s-3)(s+2K)}{(s-4)(s+2K+1)}.
\end{align*}
Hence, we obtain (3).
\par (4) Let $K \to \infty$ in (3). Then we have
\[
\lim_{K \to \infty} \frac{s+2K-2}{s+2K-1}=1.
\]
Thus we obtain
\[
Z_{2}^{\infty}(s) = \frac{s-3}{s-4}.
\]
\end{proof}

\begin{proof}[Proof of Theorem 3]
\par (1) For an integer $N \geq 1$ we prove
\[
Z_{N}^{K}(5-KN-s)=Z_{N}^{K}(s)
\]
more generally. Put $K=0$. Then we have
\begin{align*}
Z_{N}^{0}(s) &= Z(s) \\
&= \frac{(s-2)(s-3)}{(s-1)(s-4)}.
\end{align*}
This gives the functional equation
\begin{eqnarray}
Z(5-s) &= \frac{(3-s)(2-s)}{(4-s)(1-s)} \nonumber\\
&=Z(s). 
\label{2.1}
\end{eqnarray}
For an integer $K \geq 0$ we have
\[
Z_{N}^{K}(s) = \prod_{k=0}^{K}Z(s+kN).
\]
This gives
\begin{align*}
Z_{N}^{K}(5-KN-s) &= \prod_{k=0}^{K}Z(5-KN-s+kN)\\
&=\prod_{k=0}^{K}Z(5-(K-k)N-s)\\
&=\prod_{k=0}^{K}Z(5-(kN+s)).
\end{align*}
From (\ref{2.1}) we obtain
\begin{align*}
Z_{N}^{K}(5-KN-s) &= \prod_{k=0}^{K}Z(s+kN)\\
&=Z_{N}^{K}(s).
\end{align*}
\par (2) By 
\[
Z_{1}^{\infty}(s) = \frac{s-2}{s-4}
\] we have
\begin{align*}
Z_{1}^{\infty}(6-s) &= \frac{4-s}{2-s}\\
&=\frac{s-4}{s-2}\\
&= Z_{1}^{\infty}(s)^{-1}.
\end{align*}
\par (3) By the proof of (1) we obtain (3).
\par (4) By
\[
Z_{2}^{\infty}(s) = \frac{s-3}{s-4}
\] we have
\begin{align*}
Z_{2}^{\infty}(7-s)&=\frac{4-s}{3-s}\\
&=\frac{s-4}{s-3}\\
&=Z_{2}^{\infty}(s)^{-1}.
\end{align*}
\end{proof}

\begin{proof}[Proof of Theorem 4]
Since
\begin{align*}
Z(s+kN) &= (s-1+kN)^{-1}(s-2+kN)(s-3+kN)(s-4+kN)^{-1}\\
&= (k+\frac{s-1}{N})^{-1} (k+\frac{s-2}{N})(k+\frac{s-3}{N})(k+\frac{s-4}{N})^{-1}
\end{align*}
and
\[
\Gamma(x+1)=x\Gamma(x),
\]
we have
\begin{align*}
Z(s+kN)= \frac{\Gamma(\frac{s-1}{N}+k)}{\Gamma(\frac{s-1}{N}+k+1)} \frac{\Gamma(\frac{s-2}{N}+k+1)}{\Gamma(\frac{s-2}{N}+k)} \frac{\Gamma(\frac{s-3}{N}+k+1)}{\Gamma(\frac{s-3}{N}+k)} \frac{\Gamma(\frac{s-4}{N}+k)}{\Gamma(\frac{s-4}{N}+k+1)}.
\end{align*}
Hence, we obtain
\begin{align*}
&Z_{N}^{K}(s) = \prod_{k=0}^{K} Z(s+kN)\\
&=\frac{\Gamma(\frac{s-1}{N})}{\Gamma(\frac{s-1}{N}+K+1)} \frac{\Gamma(\frac{s-2}{N}+K+1)}{\Gamma(\frac{s-2}{N})} \frac{\Gamma(\frac{s-3}{N}+K+1)}{\Gamma(\frac{s-3}{N})} \frac{\Gamma(\frac{s-4}{N})}{\Gamma(\frac{s-4}{N}+K+1)}\\
&=\frac{\Gamma(\frac{s-1}{N})\Gamma(\frac{s-4}{N})}{\Gamma(\frac{s-2}{N})\Gamma(\frac{s-3}{N})} \frac{\Gamma(\frac{s-2}{N}+K+1)\Gamma(\frac{s-3}{N}+K+1)}{\Gamma(\frac{s-1}{N}+K+1)\Gamma(\frac{s-4}{N}+K+1)}.
\end{align*}
As $K \to \infty$ for $\alpha \in {\mathbb C}$ by the Stirling's formula we have
\begin{align*}
\Gamma(\alpha+K+1) &\sim \sqrt{2\pi} (\alpha + K)^{\alpha+K+\frac{1}{2}} e^{-(\alpha+K)}\\
&= \sqrt{2 \pi} K^{\alpha+K+\frac{1}{2}} (1+\frac{\alpha}{K})^{\alpha+K+\frac{1}{2}} e^{-(\alpha+K)}\\
&\sim \sqrt{2 \pi} K^{\alpha + K + \frac{1}{2}} e^{\alpha} e^{-(\alpha+K)}\\
&=\sqrt{2 \pi} K^{\alpha + K + \frac{1}{2}} e^{-K}.
\end{align*}
Using this formula to $\alpha = \frac{s-1}{N}, \frac{s-2}{N}, \frac{s-3}{N}, \frac{s-4}{N}$, we have
\[
\lim_{K \infty} \frac{\Gamma(\frac{s-2}{N}+K+1)\Gamma(\frac{s-3}{N}+K+1)}{\Gamma(\frac{s-1}{N}+K+1)\Gamma(\frac{s-4}{N}+K+1)} =1.
\]
Thus we obtain
\begin{align*}
Z_{N}^{\infty}(s) &= \lim_{K \to \infty} Z_{N}^{K}(s)\\
&=\frac{\Gamma(\frac{s-1}{N})\Gamma(\frac{s-4}{N})}{\Gamma(\frac{s-2}{N})\Gamma(\frac{s-3}{N})}.
\end{align*}
\end{proof}

\begin{proof}[Proof of Theorem 5]
By (\ref{1.1}) and (\ref{1.2}) $Z_{1}^{\infty}(s)$ and $Z_{2}^{\infty}(s)$ are rational functions. Now we assume $N \geq 3$. Then $Z_{N}^{\infty}(s)$ has poles at
\[
s=1-nN \quad (n=0,1,2,...):
\]
$\Gamma(\frac{s-1}{N})$ has poles of order $1$. $\Gamma(\frac{s-2}{N})\Gamma(\frac{s-3}{N}) \biggl|_{s=1-nN}=\Gamma(-(n+\frac{1}{N}))\Gamma(-(n+\frac{2}{N})) (\neq 0)$ is a finite value. $\Gamma(\frac{s-4}{N})\left\{
\begin{array}{cc}
\rm{has \; poles \; of \; order} \; 1 & (N=3),\\
(\neq 0) \; \rm{is \; a \; finite \; value} & (N \geq 4).
\end{array}
\right.$ Hence, $Z_{N}^{\infty}(s)$ is not a rational function when $N \geq 3$.

\end{proof}

\section{Proof of Theorem 6--8.}
\begin{proof}[Proof of Theorem 6]
\par (1) We prove Theorem 6 by the same way of the proof of Theorem 4. Since
\begin{align*}
Z(s+kN) &= (s-a+kN)^{-1}(s-b+kN)^{-1}(s-c+kN)(s-d+kN)\\
&= (k+\frac{s-a}{N})^{-1} (k+\frac{s-b}{N})^{-1} (k+\frac{s-c}{N})(k+\frac{s-d}{N})
\end{align*}
we have
\begin{align*}
Z(s+kN)= \frac{\Gamma(\frac{s-a}{N}+k)}{\Gamma(\frac{s-a}{N}+k+1)} \frac{\Gamma(\frac{s-b}{N}+k)}{\Gamma(\frac{s-b}{N}+k+1)} \frac{\Gamma(\frac{s-c}{N}+k+1)}{\Gamma(\frac{s-c}{N}+k)} \frac{\Gamma(\frac{s-d}{N}+k+1)}{\Gamma(\frac{s-d}{N}+k)}.
\end{align*}
Hence, we obtain
\begin{align*}
Z_{N}^{K}(s) &= \frac{\Gamma(\frac{s-a}{N})}{\Gamma(\frac{s-a}{N}+K+1)} \frac{\Gamma(\frac{s-b}{N})}{\Gamma(\frac{s-b}{N}+K+1)} \frac{\Gamma(\frac{s-c}{N}+K+1)}{\Gamma(\frac{s-c}{N})} \frac{\Gamma(\frac{s-d}{N}+K+1)}{\Gamma(\frac{s-d}{N})}\\
&=\frac{\Gamma(\frac{s-a}{N})\Gamma(\frac{s-b}{N})}{\Gamma(\frac{s-c}{N})\Gamma(\frac{s-d}{N})} \frac{\Gamma(\frac{s-c}{N}+K+1)\Gamma(\frac{s-d}{N}+K+1)}{\Gamma(\frac{s-a}{N}+K+1)\Gamma(\frac{s-b}{N}+K+1)}.
\end{align*}
Here as $K \to \infty$ by the Stirling's formula we have 
\begin{align*}
\Gamma(\frac{s-a}{N}+K+1) \sim \sqrt{2 \pi} K^{\frac{s-a}{N}+K+\frac{1}{2}}e^{-K},\\
\Gamma(\frac{s-b}{N}+K+1) \sim \sqrt{2 \pi} K^{\frac{s-b}{N}+K+\frac{1}{2}}e^{-K},\\
\Gamma(\frac{s-c}{N}+K+1) \sim \sqrt{2 \pi} K^{\frac{s-c}{N}+K+\frac{1}{2}}e^{-K},\\
\Gamma(\frac{s-d}{N}+K+1) \sim \sqrt{2 \pi} K^{\frac{s-d}{N}+K+\frac{1}{2}}e^{-K}.
\end{align*}
These formula give
\[
\frac{\Gamma(\frac{s-c}{N}+K+1)\Gamma(\frac{s-d}{N}+K+1)}{\Gamma(\frac{s-a}{N}+K+1)\Gamma(\frac{s-b}{N}+K+1)} \sim K^{\frac{a+b-c-d}{N}}.
\]
Hence, the convergence of $Z_{N}^{\infty}(s)$ is equivalent to $a+b=c+d$. Then we have
\[
Z_{N}^{\infty}(s) = \frac{\Gamma(\frac{s-a}{N})\Gamma(\frac{s-b}{N})}{\Gamma(\frac{s-c}{N})\Gamma(\frac{s-d}{N})}.
\]
\end{proof}

\begin{proof}[Proof of Theorem 7]
From Theorem 5 the convergence of $Z_{N}^{\infty}(s)$ is equivalent to $a+b=c+d$. So we show (2) and (3) are equivalent to $a+b=c+d$ respectively.
\par (2) 
\begin{align*}
&Z(a+b-s)=Z(s)\\
&\Leftrightarrow \frac{(a+b-c-s)(a+b-d-s)}{(a-s)(b-s)}=\frac{(s-c)(s-d)}{(s-a)(s-b)}\\
&\Leftrightarrow \frac{(s-(a+b-c))(s-(a+b-d))}{(s-a)(s-b)}=\frac{(s-c)(s-d)}{(s-a)(s-b)}\\
&\Leftrightarrow (s-(a+b-c))(s-(a+b-d)) = (s-c)(s-d)\\
&\Leftrightarrow s^2 - (2a + 2b -c -d)s+(a+b-c)(a+b-d)=s^2 -(c+d)s+cd\\
&\Leftrightarrow a+b = c+d.
\end{align*}
\par (3)
\begin{align*}
f(\frac{1}{x}) &= x^{-(a+b)}f(x)\\
&\Leftrightarrow x^{a+b} (x^{-a}+x^{-b}-x^{-c}-x^{-d})=x^{a}+x^{b}-x^{c}-x^{d}\\
&\Leftrightarrow x^{a+b-c}+x^{a+b-d}=x^{c}+x^{d}\\
&\Leftrightarrow a+b=c+d.
\end{align*}
\end{proof}

\begin{proof}[Proof of Theorem 8] 
For $f(x)=x^a+x^b-x^c-x^d$ we put 
\[
f_{N}^{\infty}(x)=\frac{f(x)}{1-x^{-N}}.
\]
Then $f_{N}^{\infty}(x)$ satisfies the absolute automorphy:
\[
f_{N}^{\infty}(\frac{1}{x}) = -x^{-(a+b+N)}f_{N}^{\infty}(x).
\]
\begin{remark} We call the absolute automorphy (see \cite{KO,KT1,KT2,T} for detail) of the condition
\[
f(\frac{1}{x}) = C x^{-D}f(x)
\]
with $C=\pm 1$ and $D \in {\mathbb R}$. Moreover we define the absolute zeta function $\zeta_{f}(s)$ of $f$ by
$$
\zeta_{f}(s) := \exp \biggl( \frac{\partial}{\partial w} Z_{f}(w,s) \biggl|_{w=0} \biggl)
$$
with
$$
Z_{f}(w,s) := \frac{1}{\Gamma(w)} \int_{1}^{\infty} f(x) x^{-s-1} (\log x)^{w-1} dx.
$$
Absolute zeta functions were studied by Soul{\'e} \cite{S} and Connes and Consani \cite{CC1}.
\end{remark}
Now we have
\begin{align*}
Z_{N}^{\infty}(s) &= \zeta_{f_{N}^{\infty}}(s)\\
&= \frac{\Gamma(\frac{s-a}{N})\Gamma(\frac{s-b}{N})}{\Gamma(\frac{s-c}{N})\Gamma(\frac{s-d}{N})}.
\end{align*}
Hence, $Z_{N}^{\infty}(s)$ is a rational function $\Leftrightarrow$ $f_{N}^{\infty}(x) \in {\mathbb Z}[x,x^{-1}]$ $\Leftrightarrow$ $f(\zeta_{N})=0$ with $\zeta_{N}=\exp(\frac{2\pi i}{N})$. Thus $Z_{N}^{\infty}(s)$ is a rational function if and only if $\zeta_{N}^{a}+\zeta_{N}^{b}=\zeta_{N}^{c}+\zeta_{N}^{d}$.
\begin{remark}
Since
\begin{align*}
\zeta_{N}^{a}+\zeta_{N}^{b} &= \zeta_{N}^{\frac{a+b}{2}} (\zeta_{N}^{\frac{a-b}{2}}+\zeta_{N}^{\frac{b-a}{2}})\\
&=2\zeta_{N}^{\frac{a+b}{2}} \cos(\frac{|a-b|}{N}\pi),\\
\zeta_{N}^{c}+\zeta_{N}^{d} &= 2\zeta_{N}^{\frac{c+d}{2}} \cos(\frac{|c-d|}{N}\pi),
\end{align*}
the condition of $N$ is given by
\begin{align*}
&\cos(\frac{|a-b|}{N}\pi)=\cos(\frac{|c-d|}{N}\pi)\\
&\Leftrightarrow \sin(\frac{|a-b|+|c-d|}{2N}\pi) \sin(\frac{|a-b|-|c-d|}{2N})=0\\
&\Leftrightarrow N|\frac{|a-b|+|c-d|}{2} or N|\frac{|a-b|-|c-d|}{2}.
\end{align*}
\end{remark}
\end{proof}

\section{Proof of Theorem 9--10.}

\begin{proof}[Proof of Theorem 9] Let $r=1$. Since
\begin{align*}
Z(s) &= \zeta_{{\mathbb G}_{m}/{\mathbb F}_{1}}(s)\\
 &= \frac{s}{s-1},
\end{align*}
we have
\begin{align*}
Z_{N}^{K}(s) &= \prod_{k=0}^{K}Z(s+kN) \\
&=\prod_{k=0}^{K}\frac{s+kN}{s-1+kN}.
\end{align*}
Hence, when $N=1$ we obtain 
\begin{align*}
Z_{1}^{K}(s) &=\prod_{k=0}^{K}\frac{s+k}{s-1+k}\\
&=\frac{s+K}{s-1}
\end{align*}
and
\begin{align*}
Z_{1}^{\infty}(s) &= \lim_{K \to \infty} Z_{1}^{K}(s)\\
&= \infty.
\end{align*}
Now we assume $N \geq 2$. We notice that
\[
Z_{N}^{k}(s) = \zeta_{f_{1,N}^{K}}(s)
\]
with
\begin{align*}
f_{1,N}^{K}(x) &= (x-1)(1+x^{-N}+\cdots+x^{-KN})\\
&=(x-1)\frac{1-x^{-(K+1)N}}{1-x^{-N}}\\
&=\frac{x+x^{-(K+1)N}-1-x^{1-(K+1)N}}{1-x^{-N}}.
\end{align*}
Using the Stirling's formula and
\begin{align*}
Z_{N}^{K}(s) &= \frac{\Gamma_{1}(s-1,(N))\Gamma_{1}(s+(K+1)N,(N))}{\Gamma_{1}(s,(N))\Gamma_{1}(s+(K+1)N-1,(N))}\\
&= \frac{\Gamma(\frac{s-1}{N})\Gamma(\frac{s}{N}+K+1)}{\Gamma(\frac{s}{N})\Gamma(\frac{s-1}{N}+K+1)},
\end{align*}
we have
\begin{align*}
Z_{N}^{\infty}(s) &= \lim_{K \to \infty} Z_{N}^{K}(s)\\
&=\infty.
\end{align*}
Next let $r \geq 2$. Then
\begin{align*}
Z(s) &= \zeta_{{\mathbb G}_{m}^{r}/{\mathbb F}_{1}}(s) \\
&=\prod_{\ell=0}^{r}(s-\ell)^{(-1)^{r-\ell+1} {\scriptsize \left(
\begin{array}{c}
r\\
\ell
\end{array}
\right)}},\\
Z_{N}^{K}(s) &=\prod_{k=0}^{K} Z(s+kN).
\end{align*}
Since 
\begin{align*}
Z(s+kN) &= \prod_{\ell=0}^{r} (s-\ell+kN)^{(-1)^{r-\ell+1}{\scriptsize \left(
\begin{array}{c}
r\\
\ell
\end{array}
\right)}}\\
&= \prod_{\ell=0}^{r} (k+\frac{s-\ell}{N})^{(-1)^{r-\ell+1}{\scriptsize \left(
\begin{array}{c}
r\\
\ell
\end{array}
\right)}}\\
&=\prod_{\ell=0}^{r}(\frac{\Gamma(\frac{s-\ell}{N}+k+1)}{\Gamma(\frac{s-\ell}{N}+k)})^{(-1)^{r-\ell+1}{\scriptsize \left(
\begin{array}{c}
r\\
\ell
\end{array}
\right)}},
\end{align*}
\begin{align*}
Z_{N}^{K}(s) &= \prod_{\ell=0}^{r} \biggl( \prod_{k=0}^{K} \frac{\Gamma(\frac{s-\ell}{N}+k+1)}{\Gamma(\frac{s-\ell}{N}+k)} \biggl)^{(-1)^{r-\ell+1}{\scriptsize \left(
\begin{array}{c}
r\\
\ell
\end{array}
\right)}}\\
&= \prod_{\ell=0}^{r} \biggl( \frac{\Gamma(\frac{s-\ell}{N}+K+1)}{\Gamma(\frac{s-\ell}{N})} \biggl)^{(-1)^{r-\ell+1}{\scriptsize \left(
\begin{array}{c}
r\\
\ell
\end{array}
\right)}}\\
&=\prod_{\ell=0}^{r} \Gamma(\frac{s-\ell}{N})^{(-1)^{r-\ell}{\scriptsize \left(
\begin{array}{c}
r\\
\ell
\end{array}
\right)}}\\
&\times \prod_{\ell=0}^{r} \Gamma(\frac{s-\ell}{N}+K+1)^{(-1)^{r-\ell+1}{\scriptsize \left(
\begin{array}{c}
r\\
\ell
\end{array}
\right)}}.
\end{align*}
As $K \to \infty$ by the Stirling's formula we have
\[
\Gamma(\frac{s-\ell}{N}+K+1) \sim \sqrt{2\pi} K^{\frac{s-\ell}{N}+K+\frac{1}{2}}e^{-K}.
\]
Hence, we have
\[
\lim_{K \to \infty} \prod_{\ell=0}^{r} \Gamma(\frac{s-\ell}{N}+K+1)^{(-1)^{r-\ell+1}{\scriptsize \left(
\begin{array}{c}
r\\
\ell
\end{array}
\right)}}=1,
\]
where we used
\[
\left\{
\begin{array}{c}
\sum_{\ell=0}^{r}(-1)^{r-\ell} \left(
\begin{array}{c}
r\\
\ell
\end{array}
\right)=0,\\
\sum_{\ell=0}^{r}(-1)^{r-\ell} \ell \left(
\begin{array}{c}
r\\
\ell
\end{array}
\right)=0.
\end{array}
\right.
\] 
Then we have
\begin{align*}
Z_{N}^{\infty}(s) &= \prod_{k=0}^{\infty} \zeta_{{\mathbb G}_{m}^{r}/{\mathbb F}_{1}}(s+kN) \\
&= \prod_{\ell=0}^{r} \Gamma(\frac{s-\ell}{N})^{(-1)^{r-\ell} {\scriptsize \left(
\begin{array}{c}
r\\
\ell
\end{array}
\right)}}
\end{align*}
\end{proof}

\begin{proof}[Proof of Theorem 10] Let 
\[
f_{r,N}^{\infty}(x)=\frac{(x-1)^{r}}{1-x^{-N}}.
\]
Since 
\begin{align*}
Z_{N}^{\infty}(s) &= \prod_{k=0}^{\infty} \zeta_{{\mathbb G}_{m}^{r}/{\mathbb F}_{1}}(s+kN)\\
&=\zeta_{f_{r,N}^{\infty}}(s),
\end{align*}
$Z_{N}^{\infty}(s)$ is a rational function $\Leftrightarrow$ $\frac{(x-1)^{r}}{1-x^{-N}} \in {\mathbb Z}[x,x^{-1}]$ $\Leftrightarrow$ $N=1$. Especially, when $N=1$ using 
\begin{align*}
f_{r,1}^{\infty}(x) &= \frac{(x-1)^{r}}{1-x^{-1}}\\
&=x(x-1)^{r-1},
\end{align*}
\[
\prod_{k=0}^{\infty}\zeta_{{\mathbb G}_{m}^{r}/{\mathbb F}_{1}}(s+k)=\zeta_{{\mathbb G}_{m}^{r-1}/{\mathbb F}_{1}}(s-1).
\]
\end{proof}

\end{document}